\documentclass[ejs,a4paper]{imsart}

\usepackage[T1]{fontenc}
\usepackage{lmodern}

\RequirePackage{amsthm,amsmath}
\RequirePackage[numbers]{natbib}
\RequirePackage[colorlinks,citecolor=blue,urlcolor=blue]{hyperref}

\usepackage[utf8]{inputenc} 
\usepackage[T1]{fontenc}    
\usepackage{hyperref}       
\usepackage{url}            
\usepackage{booktabs}       
\usepackage{amsfonts}       
\usepackage{nicefrac}       
\usepackage{microtype}      

\usepackage[utf8]{inputenc}
\usepackage{a4wide}
\usepackage[most]{tcolorbox}
\usepackage{amssymb}
\usepackage{amsfonts,mathrsfs,amsmath}
\usepackage{pifont,dsfont,booktabs}
\usepackage{breakcites,microtype}

\usepackage[export]{adjustbox}
\usepackage{graphicx}

\bibpunct{(}{)}{;}{a}{,}{,}
\usepackage{hyperref}
\hypersetup{
  colorlinks = true,
  urlcolor = blue,
  linkcolor = blue,
  citecolor = blue,
}

\usepackage{cleveref}
\usepackage{autonum}
\usepackage{chngcntr}

\usepackage{enumitem}
\usepackage{nicefrac}

\usepackage{caption}
\crefname{algocf}{alg.}{algs.}
\Crefname{algocf}{Algorithm}{Algorithms}
 
\def\bs{\boldsymbol}
\def\dTV{d_{\rm TV}}
\def\dKL{d_{\rm KL}}

\def\RR{\mathbb R}

\def\bx{\boldsymbol x}
\def\by{\boldsymbol y}

\def\btheta{\boldsymbol\theta}

\def\bY{\boldsymbol Y\!}
\def\bX{\boldsymbol X}
\def\bZ{\boldsymbol Z}

\def\bQ{\boldsymbol Q}
\def\bP{\boldsymbol P}
\def\bPn{{\boldsymbol P}\!_n}

\def\bfE{\mathbf E}
\def\bfP{\mathbf P}

\def\bu{\boldsymbol u}
\def\bv{\boldsymbol v}

\def\be{\boldsymbol e}

\def\fcar{\mathds{1}}

\def\var{\mathbf{Var}}

\def\calM{\mathcal M}

\def\simiid{\overset{\text{iid}}{\sim}}

\def\hat{\widehat}

\def\epsilon{\varepsilon}

\usepackage{pgfpages}
\hypersetup{pdfpagemode=UseNone}
\usepackage[most]{tcolorbox}
\tcbuselibrary{skins,breakable,xparse}

\NewTColorBox[auto counter]{thm}{+O{}}{ %
enhanced,
drop fuzzy shadow,
colframe=green!20!black,
colback=yellow!10!white,
coltitle=green!10!black,
fonttitle=\bfseries,
title={Theorem~\thetcbcounter},
label={theorem@\thetcbcounter},
attach title to upper=\quad,
lowerbox=ignored,
#1
}

\NewTotalTColorBox{\tproof}{mm}{ %
boxsep=0mm,
left=0mm,
right=0mm,
enhanced,breakable,
colframe=white,
colback=white,
coltitle=black,
fonttitle=\bfseries,
title={Proof of Theorem~\ref{theorem@#1} on page~\pageref{theorem@#1}},
label={proof@#1},
attach title to upper=\par,
}{\input{#2}}

\NewTColorBox[auto counter]{prop}{+O{}}{ %
enhanced,
drop fuzzy shadow,
colframe=green!20!black,
colback=yellow!10!white,
coltitle=green!10!black,
fonttitle=\bfseries,
title={Proposition~\thetcbcounter},
label={prop@\thetcbcounter},
attach title to upper=\quad,
lowerbox=ignored,
#1
}

\NewTotalTColorBox{\pproof}{mm}{ %
boxsep=0mm,
left=0mm,
right=0mm,
enhanced,breakable,
colframe=white,
colback=white,
coltitle=black,
fonttitle=\bfseries,
title={},
phantomlabel={proof@#1},
attach title to upper=\par,
}{\textbf{Proof of Proposition~\ref{prop@#1} on page~\pageref{prop@#1} }\input{#2}}

\tcbset{no proof/.style={no recording,after upper=}}

\begin{document}
\begin{frontmatter}
\title{Confidence regions and minimax rates in outlier-robust estimation on the probability simplex}
\runtitle{Outlier-robust estimation on the probability simplex}

\begin{aug}
\author{\fnms{Amir-Hossein} \snm{Bateni}\ead[label=e1]{amirhossein.bateni@ensae.fr}}
\and
\author{\fnms{Arnak~S.} \snm{Dalalyan}\ead[label=e2]{arnak.dalalyan@ensae.fr}}

\address{5 Avenue Le Chatelier, 91120 Palaiseau, France\\
}

\runauthor{A. Bateni and A. Dalalyan}

\affiliation{ENSAE-CREST}

\end{aug}


\begin{abstract}
We consider the problem of estimating the mean of a 
distribution supported by the $k$-dimensional probability
simplex in the setting where an $\varepsilon$ fraction of
observations are subject to adversarial corruption. A 
simple particular example is the problem of estimating the
distribution of a discrete random variable. Assuming 
that the discrete variable takes $k$ values, the unknown
parameter $\bs \theta$ is a $k$-dimensional vector belonging 
to the probability simplex. We first describe various settings 
of contamination and discuss the relation between these settings. 
We then establish minimax rates when the quality of estimation 
is measured by the total-variation distance, the Hellinger 
distance, or the $\mathbb L^2$-distance between two probability 
measures. We also provide confidence regions for the unknown 
mean that shrink at the minimax rate. Our analysis reveals 
that the minimax rates associated to these three distances 
are all different, but they are all attained by the sample 
average. Furthermore, we show that the latter is adaptive 
to the possible sparsity of the unknown vector. Some 
numerical experiments illustrating our theoretical findings 
are reported.
\end{abstract}

\setcounter{tocdepth}{1}
\tableofcontents
\end{frontmatter}

\parskip=5pt
\section{Introduction}\label{sec:intro}

Assume $\bX_1,\ldots,\bX_n$ are $n$ independent random variables
taking their values in the $k$-dimensional probability simplex 
$\Delta^{k-1} = \{\bv\in\RR_+^k:v_1+\ldots+v_k=1\}$. Our goal is to
estimate the unknown vector $\bs \theta = \bfE[\bX_i]$ in the case
where the observations are contaminated by outliers. In this
introduction, to convey the main messages, we limit ourselves to 
the Huber contamination model, although our results apply to 
the more general adversarial contamination. Huber's contamination 
model assumes that there are two probability measures $\bs P$, $\bs Q$ 
on $\Delta^{k-1}$ and a real $\varepsilon\in[0,1/2)$ such that 
$\bX_i$ is drawn from 
\begin{align}
\bs P_i = (1-\varepsilon)\bs P + \varepsilon \bs Q,\qquad \forall 
i\in \{1,\ldots,n\}.
\end{align}
This amounts to assuming that $(1-\epsilon)$-fraction of observations,
called inliers, are drawn from a reference distribution $\bs P$, whereas
$\epsilon$-fraction of observations are outliers and are drawn from 
another distribution $\bs Q$. In general, all the three parameters 
$\bs P$, $\bs Q$ and $\varepsilon$ are unknown. The parameter of
interest is some functional (such  as the mean, the standard deviation,
etc.) of the reference distribution $\bs P$, whereas $\bs Q$ and
$\varepsilon$ play the role of nuisance parameters.

When the unknown parameter lives on the probability simplex, there
are many appealing ways of defining the risk. We focus on the 
following three metrics: total-variation, Hellinger and $\mathbb L^2$ 
distances
\footnote{We write
$\|\bu\|_q = (\sum_{j=1}^k |u_j|^q)^{1/q}$ and $\bu^q = (u_1^q,\ldots,u_k^q)$ for any $\bu\in\mathbb R_+^k$ and $q> 0$.}
\begin{align}
d_{\text{TV}}(\hat{\bs \theta},\bs \theta) :=\nicefrac1{2} \|\hat{\bs \theta} - \bs \theta\|_{1}
,\quad
d_{\text{H}}(\hat{\bs \theta},\bs \theta) := 
\|\hat{\bs \theta}{}^{1/2} - \bs \theta^{1/2}\|_{2}
,\quad
d_{\mathbb L^2}(\hat{\bs \theta},\bs \theta) := \|\hat{\bs \theta} - \bs \theta\|_{2} 
.
\end{align}
The Hellinger distance above is well defined when the estimator
$\hat{\bs \theta}$ is non-negative, which will be the case
throughout this work. We will further assume that the dimension
$k$ may be large, but the vector $\bs\theta$ is $s$-sparse, for
some $s\le k$, \textit{i.e.}  $\#\{j:\theta_j\neq 0\}\le s$. 
Our main interest is in constructing confidence regions and evaluating the minimax risk
\begin{align}\label{mmx_def}
\mathfrak{R}_\square(n,k,s,\varepsilon) := 
\inf_{\bar{\bs \theta}_n}\sup_{\bs P,\bs Q}
\bfE[d_\square(\bar{\bs \theta}_n,\bs \theta)],
\end{align}
where the \emph{inf} is over all estimators $\bar{\bs\theta}_n$
built upon the observations $\bX_1,\ldots,\bX_n\simiid
(1-\varepsilon)\bs P+ \varepsilon \bs Q$ and the \emph{sup} 
is over all distributions $\bs P$, $\bs Q$ on the probability 
simplex such that the mean $\bs\theta$ of ${\bs P}$ is $s$-sparse. 
The subscript $\square$ of $\mathfrak{R}$ above refers to the 
distance used in the risk, so that $\square$ is TV, H, or 
$\mathbb L^2$. 

The problem described above arises in many practical
situations. One example is an election poll: each participant
expresses his intention to vote for one of $k$ candidates.
Thus, each $\theta_j$ is the true proportion of electors of
candidate $j$. The results of the poll contain outliers,
since some participants of the poll prefer to hide their true
opinion. Another example, still related to elections, is
the problem of counting votes across all constituencies. 
Each constituency communicates a vector of proportions to
a central office, which is in charge of computing the 
overall proportions. However, in some constituencies
(hopefully a small fraction only) the results are rigged.  Therefore, the set of observed vectors contains some outliers.

We intend to provide non-asymptotic upper and 
lower bounds on the minimax risk that match up to numerical constants. In addition, we will provide confidence regions of the form $B_\square(\hat{\bs\theta}_n, r_{n,\varepsilon,\delta})= \{\bs\theta : d_\square(\hat{\bs\theta}_n,\bs\theta)\le r_{n,\varepsilon,\delta}\}$ containing the true parameter with probability at least $1-\delta$ and such that the radius $r_{n,\varepsilon,\delta}$ goes to zero at the same rate as the corresponding minimax risk.

When there is no outlier, \textit{i.e.}, $\varepsilon =0$, 
it is well known that the sample mean
\begin{align}
\bar{\bs X}_n := \frac1n \sum_{i=1}^n \bX_i
\end{align}
is minimax-rate-optimal and the rates corresponding to various 
distances are
\begin{align}
\mathfrak{R}_{\mathbb L^2}(n,k,s,0) \asymp (1/n)^{1/2} \quad
\text{and}\quad
\mathfrak{R}_{\square}(n,k,s,0) \asymp (s/n)^{1/2}\quad
\text{for } \square\in\{\text{TV},\text{H}\}.
\end{align}
This raises several questions in the setting where data 
contains outliers. In particular, the following three questions
will be answered in this work:
\begin{description}
\item[Q1.] How the risks $\mathfrak{R}_\square$ depend on 
$\varepsilon$? What is the largest proportion of outliers for which the minimax rate is the same
as in the outlier-free case ?
\item[Q2.] Does the sample mean remain optimal in the 
contaminated setting?
\item[Q3.] What happens if the unknown parameter $\bs\theta$ 
is $s$-sparse ?
\end{description}
The most important step for answering these questions is to
show that 
\begin{align}
\mathfrak{R}_{\text{TV}}(n,k,s,\varepsilon)&\asymp (s/n)^{1/2}+\varepsilon,\\ 
\mathfrak{R}_{\text{H}}(n,k,s,\varepsilon)&\asymp (s/n)^{1/2}+\varepsilon^{1/2},\\ 
\mathfrak{R}_{\mathbb L^2}(n,k,s,\varepsilon)&\asymp (1/n)^{1/2}+\varepsilon.
\end{align}
It is surprising to see that all the three rates are
different leading to important discrepancies in the answers
to the second part of question Q1 for different distances.
Indeed, it turns out that the minimax rate is not
deteriorated if the proportion of the outliers is smaller
than $(s/n)^{1/2}$ for the TV-distance, $s/n$ for the
Hellinger distance and $(1/n)^{1/2}$ for the $\mathbb L^2$
distance. Furthermore, we prove that the sample mean is
minimax rate optimal. Thus, even when the proportion of
outliers $\varepsilon$ and the sparsity $s$ are known, it 
is not possible to improve upon the sample mean. In 
addition, we show that all these claims hold true for the
adversarial contamination and we provide corresponding confidence regions.

The rest of the paper is organized as follows. \Cref{sec:2} 
introduces different possible ways of modeling data sets contaminated  by outliers. Pointers to relevant prior work
are given in \Cref{sec:3}. Main theoretical results and their
numerical illustration are reported in \Cref{sec:4} and
\Cref{sec:5}, respectively. \Cref{sec:6} contains a brief
summary of the obtained results and their consequences,
whereas the proofs are postponed to the appendix.

\section{Various models of contamination}\label{sec:2}

Different mathematical frameworks have been used in the
literature to model the outliers. We present here five of
them, from the most restrictive one to the most general, 
and describe their relationship. We present these frameworks
in the general setting when the goal is to estimate the
parameter $\btheta^*$ of a reference distribution $\bP_{\btheta^*}$ when $\varepsilon$ proportion of the
observations are outliers.

\subsection{Huber's contamination}

The most popular framework for studying robust estimation
methods is perhaps the one of Huber's contamination. In this
framework, there is a distribution $\bQ$ defined on the same
space as the reference distribution $\bP_{\btheta^*}$ such
that all the observations $\bX_1,\ldots,\bX_n$ are
independent and drawn from the mixture distribution
$\bP_{\varepsilon,\btheta^*,\bQ} := (1-\varepsilon)
\bP_{\btheta^*} + \varepsilon \bQ$.

This corresponds to the following mechanism: one decides with
probabilities  $(1-\varepsilon, \varepsilon)$ whether a given
observation is an inlier or an outlier. If the decision is
made in favor of being inlier, the observation is drawn from 
$\bP_{\btheta^*}$, otherwise it is drawn from $\bQ$. More
formally, if we denote by $\hat O$ the random set of
outliers, then conditionally to $\hat O = O$,
\begin{align}\label{eq:1}
\{\bX_i:i\not\in O\}\simiid \bP_{\btheta^*},\quad 
\{\bX_i:i\in O\}\simiid \bQ,\quad
\{\bX_i:i\in O\}\perp\!\!\!\perp \{\bX_i:i\not\in O\},
\end{align}
for every $O\subset \{1,\ldots,n\}$. Furthermore, for every
subset $O$ of the observations, we have $\bP(\hat O = O) =
(1-\varepsilon)^{n-|O|}\varepsilon^{|O|}$. We denote
by\footnote{The superscript HC refers to the Huber's 
contamination}  $\mathcal M_n^{\text{HC}} (\varepsilon,
\btheta^*)$ the set of joint probability distributions 
$\bPn$ of the random variables $\bX_1,\ldots,\bX_n$
satisfying the foregoing condition. 

\subsection{Huber's deterministic contamination}

The set of outliers as well as the number of outliers in
Huber's model of contamination are random. This makes it
difficult to compare this model to the others that will be
described later in this section. To cope with this, we 
define here another model, termed Huber's deterministic
contamination. As its name indicates, this new model has the
advantage of containing a deterministic number of outliers,
in the same time being equivalent to Huber's contamination 
in a sense that will be made precise below.

We say that the distribution $\bPn$ of $\bX_1,\ldots,\bX_n$
belongs to the Huber's deterministic contamination model
denoted by $\mathcal M_n^{\text{HDC}} (\varepsilon,
\btheta^*)$, if there are a set $O\subset \{1,\ldots,n\}$ of
cardinality at most $n\varepsilon$ and a distribution $\bQ$
such that \eqref{eq:1} is true. The apparent similarity of
models $\mathcal M_n^{\text{HC}}(\varepsilon,\btheta^*)$
and $\mathcal M_n^{\text{HDC}}(\varepsilon,\btheta^*)$ can
also be formalized mathematically in terms of the orders of
magnitude of minimax risks. To ease notation, we let $R_d^{\square}(n,\varepsilon,\Theta,\hat\btheta)$ to be the 
worst-case risk of an estimator $\hat\btheta$, where
$\square$ is either HC or HDC. More precisely, for $\mathcal
M_n^{\square}(\varepsilon,\Theta) := \cup_{\btheta\in
\Theta}\mathcal M_n^{\square}(\varepsilon,\btheta)$, 
we set\footnote{The subscript $d$ refers to the distance 
$d$ used in the definition of the risk.}
\begin{align}
    R_d^{\square}(n,\varepsilon,\Theta,\hat\btheta)
    :=\sup_{P_n\in\mathcal M_n^{\square}(\varepsilon,\Theta)}
    \bfE[d(\hat\btheta,\btheta^*)].
\end{align}
This definition assumes that the parameter space $\Theta$ is endowed with a 
pseudo-metric $d:\Theta\times\Theta\to \RR_+$.  
When $\Theta = \{\btheta^*\}$ is a singleton, we write
$R_{d,n}^{\square}(\varepsilon,\btheta^*,\hat\btheta)$ instead of
$R_d^{\square}(n,\varepsilon,\{\btheta^*\},\hat\btheta)$.

\tcbstartrecording

\begin{prop}\label{prop:1}
Let $\hat\btheta_n$ be an arbitrary estimator of $\btheta^*$. For any 
$\varepsilon\in(0,1/2)$,
\begin{align}
R_d^{\text{HC}}(n,\varepsilon,\btheta^*,\hat\btheta_n) \le R_{d,n}^{\text{HDC}}(2\varepsilon,
\btheta^*,\hat\btheta_n)  + e^{-n\varepsilon/3}\!R_{d,n}^{\text{HDC}}(1,\btheta^*,\hat\btheta_n),
\label{ineq:1}\\
\sup_{\bP_n\in \calM_n^{\text{HC}}(\varepsilon,\btheta^*)} 
r\bfP\big(d(\hat\btheta_n,\btheta^*)>r\big) \le R_{d,n}^{\text{HDC}}
(2\varepsilon,\btheta^*,\hat\btheta_n) + re^{-n\varepsilon/3}.
\label{ineq:1b}
\end{align}
{\par\hfill\textcolor{green!40!black} %
{\itshape Proof in the appendix, page~\pageref{proof:prop1}}}
\end{prop}

Denote by $\mathcal D_\Theta$ the diameter of  $\Theta$,
$\mathcal D_\Theta :=\max_{\btheta,\btheta'} d(\btheta,
\btheta')$. Proposition~\ref{prop@2}  implies that
\begin{align}\label{eq:3}
\inf_{\hat\btheta_n} R_d^{\text{HC}}(n,\varepsilon,\Theta,\hat\btheta_n) \le 
\inf_{\hat\btheta_n} R_d^{\text{HDC}}(n,2\varepsilon,
\Theta,\hat\btheta_n)  + e^{-n\varepsilon/3}\!\mathcal D_\Theta.
\end{align} 
When $\Theta$ is bounded, the last term is typically of
smaller order than the minimax risk over
$\calM^{\text{HDC}}_n (2\varepsilon,\Theta)$. Therefore, the 
minimax rate of estimation in Huber's model is not slower
than the minimax rate of estimation in Huber's deterministic
contamination model. This entails that a lower bound on the
minimax risk established in HC-model furnishes a lower
bound in HDC-model. 

\subsection{Oblivious contamination} 

A third model of contamination that can be of interest is 
the oblivious contamination. In this model, it is assumed
that the set $O$ of cardinality $o$ and the joint 
distribution $\bQ_O$ of outliers are determined in advance,
possibly based on the knowledge of the reference distribution
$\bP_{\btheta^*}$. Then, the outliers $\{\bX_i:i\in O\}$ are
drawn randomly from $\bQ_O$ independently of the inliers 
$\{\bX_i:i\in O^c\}$. The set of all the joint distributions
$\bP_n$ of random variables $\bX_1,\ldots,\bX_n$ generated by
such a mechanism will be denoted by $\mathcal M_n^{\text{OC}}
(\varepsilon,\btheta^*)$. The model of oblivious
contamination is strictly more general than that of Huber's 
deterministic contamination, since it does not assume that
the outliers are iid. Therefore, the minimax risk over $\mathcal M_n^{\text{OC}}(\varepsilon,\Theta)$
is larger than the minimax risk over $\mathcal M_n^{\text{HDC}}(\varepsilon,\Theta)$:
\begin{equation}
\inf_{\hat\btheta_n} R_d^{\text{HDC}}(n,\varepsilon,\Theta,\hat\btheta_n) 
\le \inf_{\hat\btheta_n} R_d^{\text{OC}} (n, \varepsilon,
\Theta,\hat\btheta_n).
\end{equation}
The last inequality holds true for any set $\Theta$, any
contamination level $\varepsilon\in(0,1)$ and any sample 
size.

\subsection{Parameter contamination}

In the three models considered above, the contamination acts on
the observations. One can also consider the case where the parameters of 
the distributions of some observations are contaminated. More precisely, 
for some set $O\subset\{1,\ldots,n\}$ selected in advance (but unobserved), the outliers $\{\bX_i:i\in O\}$ are independent and independent of the inliers $\{\bX_i:i\in O^c\}$. Furthermore, each outlier $\bX_i$ is drawn from a distribution $\bQ_i = \bP_{\btheta_i}$ belonging to the same family as the reference distribution, but corresponding to a contaminated parameter $\btheta_i\not=\btheta^*$. Thus, the joint distribution of the observations can be written as $(\bigotimes_{i\in O^c} \bP_{\btheta^*})\otimes(\bigotimes_{i\in O}\bP_{\btheta_i})$. The set of all such distributions $\bP_n$  will be denoted by 
$\mathcal M_n^{\text{PC}}(\varepsilon,\btheta^*)$, where PC refers to ``parameter contamination''.

\subsection{Adversarial contamination}

The last model of contamination we describe in this work, the adversarial contamination, is the most general one. It corresponds to the following two-stage data generation mechanism. In a first stage, iid random variables $\bY_1,\ldots,\bY_n$ are generated from a reference distribution $\bP_{\btheta^*}$. In a second stage, an adversary having access to  $\bY_1,\ldots,\bY_n$ chooses a (random) set $\hat O$ of (deterministic) cardinality $s$ and arbitrarily modifies data points $\{\bY_i:i\in \hat O\}$. The resulting sample, $\{\bX_i:i= 1,\ldots, n\}$, is revealed to the Statistician. In this model, we have 
$\bX_i = \bY_i$ for $i\not\in\hat O$. However, since $\hat O$ is random and potentially dependent of $\bY_{1:n}$, it is not true that conditionally to $\hat O = O$, $\{\bX_i:i\in O^c\}$  are iid drawn
from $\bP_{\btheta^*}$ (for any deterministic set $O$ of cardinality
$o$). 

We denote by $\mathcal M_n^{\text{AC}}(\varepsilon,\btheta^*)$ the set of all the joint distributions $\bP_n$ of all the sequences $\bX_1,\ldots,\bX_n$ generated by the aforementioned two-stage mechanism. This set $\mathcal M_n^{\text{AC}}(\varepsilon,\btheta^*)$ is larger than all the four sets of contamination introduced in this section. Therefore, the following inequalities hold:
\begin{align}
    \inf_{\hat\btheta_n} R_d^{\text{PC}}(n,\varepsilon,\Theta,\hat\btheta_n) \le 
		\inf_{\hat\btheta_n} R_d^{\rm OC}(n,\varepsilon,\Theta,\hat\btheta_n) \le 
		\inf_{\hat\btheta_n} R_d^{\text{AC}}(n,\varepsilon,\Theta,\hat\btheta_n),
\end{align}
for any $n$, $\varepsilon$, $\Theta$ and any distance $d$. 
\begin{figure}
        \def\svgwidth{0.8\textwidth}
				\centerline{\input{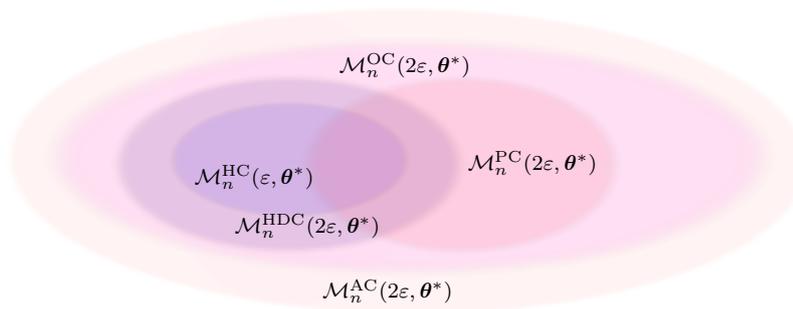}}
\caption{Visual representation of the hierarchy between various contamination model. Note 
that the inclusion of $\calM_n^{\text{HC}}(\varepsilon,\btheta^*)$ in 
$\calM_n^{\text{HDC}}(2\varepsilon,\btheta^*)$ is somewhat heuristic, based on the
relation on the worst-case risks reported in Proposition~\ref{prop@1}.}
\label{fig:1}
\end{figure}
\subsection{Minimax risk ``in expectation'' versus ``in deviation''}

Most prior work on robust estimation focused on establishing upper bounds on the minimax risk in deviation\footnote{We call a risk bound in deviation any bound on the distance $d(\hat\btheta,\btheta^*)$ that holds true with a probability close to one, for any parameter value $\btheta^*\in\Theta$.}, as opposed to the minimax risk in expectation defined by \eqref{mmx_def}. One of the reasons for dealing with the deviation is that it makes the minimax risk meaningful for models\footnote{This is the case, for instance, of the Gaussian model with Huber's contamination.} having random number of outliers and unbounded parameter space $\Theta$. The formal justification of this claim is provided by the following result.

\begin{prop}\label{prop:2}
Let $\Theta$ be a parameter space such that $\mathcal D_\Theta= \sup_{\btheta,\btheta'\in\Theta} d(\btheta,\btheta')=+\infty$. Then, for every estimator $\hat\theta_n$, every $\varepsilon>0$ and $n\in\mathbb{N}$, we have 
$R_d^{\text{HC}}(n,\varepsilon,\Theta,\hat\btheta_n) = +\infty$. 
{\par\hfill\textcolor{green!40!black} %
{\itshape Proof in the appendix, page~\pageref{proof:prop2}}}
\end{prop}

This result shows, in particular, that the last term in \eqref{eq:3}, involving the  
diameter of $\Theta$ is unavoidable. Such an explosion of the minimax risk occurs because
Huber's model allows the number of outliers to be as large as $n/2$ with a strictly positive
probability. One approach to overcome this shortcoming is to use the minimax risk in 
deviation. Another approach is to limit theoretical developments to the models HDC, PC, 
OC or AC, in which the number of outliers is deterministic. 

\section{Prior work} \label{sec:3}

Robust estimation is an area of active research in Statistics since at least five decades  \citep{huber1964,Tukey75,donoho83,donoho1992,Rousseeuw99}. Until very 
recently, theoretical guarantees were almost exclusively formulated in terms of 
the notions of  breakdown point, sensitivity curve, influence function, etc. 
These notions are 
well suited for accounting for gross outliers, observations that deviate 
significantly from the data points representative of an important fraction 
of data set. 

More recently, various authors investigated \citep{nguyen2013,Dal12,Chen} 
the behavior of the risk of robust estimators as a function of the rate 
of contamination $\varepsilon$. A general methodology for parametric models
subject to Huber's contamination was developed in \cite{CGR15, CGR16}. 
This methodology allowed for determining the rate of convergence of 
the minimax risk as a function of the sample size $n$, dimension $k$ 
and the rate of contamination $\varepsilon$. An interesting 
phenomenon was discovered: in the problem of robust estimation 
of the Gaussian mean, classic robust estimators such as the 
coordinatewise median or the geometric median do not attain the 
optimal rate $(k/n)^{1/2} + \varepsilon$. This rate is provably
attained by Tukey's median, the computation of which is costly in 
a high dimensional setting. 

In the model analyzed in this paper, we find the same minimax rate, 
$(k/n)^{1/2} + \varepsilon$, only when the total-variation distance is 
considered. A striking difference is that this rate is attained by the
sample mean which is efficiently computable in any dimension. This property
is to some extent similar to the problem of robust density estimation 
\citep{Gao2017}, in which the standard kernel estimators are minimax optimal
in contaminated setting. 


Computational intractability of Tukey's median motivated a large number
of studies that aimed at designing computationally tractable methods with nearly
optimal statistical guarantees. Many of these works went beyond
Huber's contamination by considering parameter contamination models 
\citep{Kush,Collier2017,Carpentier}, oblivious contamination \citep{Feng,LRV}
or adversarial contamination \citep{DiakonikolasKKL16,Balakrishnan,
DiakonikolasKK017,Diakonikolas0S18}. Interestingly, in the problem of
estimating the Gaussian mean, it was proven that the minimax rates under 
adversarial contamination are within a factor at most logarithmic in $n$ 
and $k$ of the minimax rates under Huber's contamination\footnote{All these
papers consider the risk in deviation, so that the minimax risk under
Huber's contamination is finite.}. While each of 
the aforementioned papers introduced clearly the conditions on the 
contamination, to our knowledge, none of them described different 
possible models  and the relationship between them. 

Another line of growing literature on robust estimation aims at robustifying 
estimators and prediction methods to heavy tailed
distributions, see \citep{audibert2011, minsker2015, Donoho2016,devroye2016,joly2017,minsker2018,lugosi2019,
lecue2017,chinot2018}. The results of those papers are of 
a different nature, as compared to the present work, not 
only in terms of the goals, but also in terms of mathematical
and algorithmic tools.  

\section{Minimax rates on the ``sparse'' simplex 
and confidence regions}\label{sec:4}

We now specialize the general setting of \Cref{sec:2} to a 
reference distribution $\bs P$, with expectation $\bs\btheta^*$,
defined on the simplex $\Delta^{k-1}$. Along with this 
reference model describing the distribution of inliers, we 
will use different models of contamination. More precisely, 
we will establish upper bounds on worst-case risks of the 
sample mean in the most general, adversarial, contamination 
setting. Then, matching lower bounds will be provided for 
minimax risks under Huber's contamination.

\subsection{Upper bounds: worst-case risk of the sample mean} 

We denote by $\Delta_s^{k-1}$ the set of all $\bv\in \Delta^{k-1}$ 
having  at most $s$ non-zero entries.

\begin{thm}\label{thm:1}
For every triple of positive integers $(k,s,n)$ and for 
every $\epsilon\in[0,1]$, the sample mean $\bar{\bs X}_n 
:= \frac1n \sum_{i=1}^n \bX_i$ satisfies
\begin{align}
    R_{\text{TV}}^{\text{AC}}(n,\epsilon,\Delta_s^{k-1},
    \bar{\bs X}_n) &\le (s/n)^{1/2} + 2\epsilon ,\\[5pt]
    R_{\text{H}}^{\text{AC}}(n,\epsilon,\Delta_s^{k-1},
    \bar{\bs X}_n) &\le
    (s/n)^{1/2} + 2\,\epsilon^{1/2} ,\\[5pt]
    R_{\mathbb{L}^2}^{\text{AC}}(n,\epsilon,\Delta_s^{k-1},
    \bar{\bs X}_n) &\le(1/n)^{1/2} + \sqrt{2}\,\epsilon.
\end{align}
{\par\hfill\textcolor{green!40!black} %
{\itshape Proof in the appendix, page~\pageref{proof:thm1}}}
\end{thm}

An unexpected and curious phenomenon unveiled by this theorem 
is that all the three rates are different. As a consequence, 
the answer to the question ``what is the largest possible number 
of outliers, $o_d^*(n,s)$, that does not impact the minimax rate 
of estimation of $\bs\theta^*$?'' crucially depends on the 
considered distance $d$. Taking into account the relation
$\varepsilon=o/n$, we get
\begin{align}
    o^*_{\text{TV}}(n,s) \asymp (ns)^{1/2},\quad
    o^*_{\text{H}}(n,s) \asymp s,\quad
    o^*_{\mathbb L^2}(n,s) \asymp n^{1/2}.
\end{align}
Furthermore, all the claims concerning the total variation
distance, in the considered 
model, yield
corresponding claims for the Wasserstein distances $W_q$, 
for every $q\ge 1$. Indeed, one can see an element 
$\bs\theta\in\Delta^{k-1}$ as the probability distribution
of a random vector $\bX$ taking values in the finite set 
$\mathcal A =\{\be_1,\ldots,\be_k\}$ of vectors of the 
canonical basis of $\mathbb R^k$. Since these vectors satisfy  
$\|\be_j-\be_{j'} \|^2_2 =2\mathds 1(j\neq j')$, we have
\begin{align} \label{WqTV}
W_q^q(\bs \theta,\bs \theta') &= \inf_{\Gamma} \bfE_{(\bX,\bX')\sim\Gamma}[\|\bX-\bX'\|_2^{q}]\\
&=
\inf_{\Gamma}  2^{q/2}\bP(\bX\neq\bX') = 2^{q/2}\|\bs \theta - \bs \theta'\|_{\text{TV}},
\end{align}
where the \emph{inf} is over all joint distributions $\Gamma$ on 
$\mathcal A\times\mathcal A$ having marginal distributions 
$\bs \theta$ and $\bs \theta'$. This implies that
\begin{align}
    R_{W_q}^{\text{AC}}(n,\epsilon,\Delta_s^{k-1}) &\le
    \sqrt{2}\big\{ (s/n)^{1/2} + 2\epsilon \big\}^{1/q},
    \qquad \forall q\ge 1.\label{Wq}
\end{align}
In addition, since the $\mathbb L_2$ norm is an upper
bound on the $\mathbb L_\infty$-norm, we have 
$R_{\mathbb L_\infty}^{\text{AC}}(n,\epsilon,\Delta^{k-1})
\le (1/n)^{1/2}+\sqrt{2}\,\varepsilon$. Thus, we have
obtained upper bounds on the risk of the sample mean for 
all commonly used distances on the space of probability
measures. 

\subsection{Lower bounds on the minimax risk}

A natural question, answered in the next theorem,
is how tight are the  upper bounds obtained in the last 
theorem. More importantly, one can wonder whether there is 
an estimator that has a worst-case risk of smaller order 
than that of the sample mean.  

\begin{thm}\label{thm:2}
There are universal constants $c>0$ and $n_0$, such that for 
any integers $k\ge 3$, $s\le k\wedge n$, $n\ge n_0$ 
and for any $\epsilon\in[0,1]$, we have
\begin{align}
    \inf_ {\bar{\bs \theta}_n}R_{\text{TV}}^{\text{HC}}(n,\epsilon,
    \Delta^{k-1}_s,\bar{\bs \theta}_n) &\ge 
    c\{(s/n)^{1/2} + \epsilon\},\\[2pt]
    \inf_ {\bar{\bs \theta}_n}R_{\text{H}}^{\text{HC}}
    (n,\epsilon,\Delta^{k-1}_s,\bar{\bs \theta}_n) 
    &\ge c\{(s/n)^{1/2} + \epsilon^{1/2}\} ,\\[2pt]
    \inf_ {\bar{\bs \theta}_n} R_{\mathbb{L}^2}^{\text{HC}}
    (n,\epsilon,\Delta^{k-1}_s,\bar{\bs \theta}_n) 
    &\ge c\{(1/n)^{1/2} + \epsilon\},
\end{align}
where $\inf_ {\bar{\bs \theta}_n}$ stands for the infimum over all measurable functions
$\bar{\bs \theta}_n$ from $(\Delta^{k-1})^{n}$ to $\Delta^{k-1}$.
{\par\hfill\textcolor{green!40!black} %
{\itshape Proof in the appendix, page~\pageref{proof:thm2}}}
\end{thm}

The main consequence of this theorem is that whatever the 
contamination model is (among those described in \Cref{sec:2}), 
the rates obtained for the MLE in \Cref{thm:1} are minimax optimal. 
Indeed, Theorem~\ref{theorem@2} yields this claim for Huber's 
contamination. For Huber's deterministic contamination and
and the TV-distance, on the one hand, we have 
\begin{align}
    R_{\text{TV}}^{\text{HDC}}(n,\epsilon,\Delta_s^{k-1},
    \bar{\bs \theta}_n) &\ge R_{\text{TV}}^{\text{HDC}}(n,0,\Delta_s^{k-1},
    \bar{\bs \theta}_n)\\
    &\stackrel{(1)}{=}
    R_{\text{TV}}^{\text{HC}}(n,0,\Delta_s^{k-1},
    \bar{\bs \theta}_n)\stackrel{(2)}{\ge} c(s/n)^{1/2},
\end{align}
where (1) uses the fact that for $\varepsilon=0$ all the sets 
$\calM_n^{\square}(\varepsilon,\btheta^*)$ are equal, while (2) 
follows from the last theorem. On the other hand, in view of Proposition~\ref{prop@1},
for $\varepsilon \ge (6/n)\log(8n/c)$ (implying that $2e^{-n\varepsilon/6}\le (c/4)\varepsilon$), 
\begin{align}
    R_{\text{TV}}^{\text{HDC}}(n,\epsilon,\Delta_s^{k-1},
    \bar{\bs \theta}_n) 
    &\ge R_{\text{TV}}^{\text{HC}}(n,\varepsilon/2,\Delta_s^{k-1},
    \bar{\bs \theta}_n) - 2e^{-n\varepsilon/6}\\
    &\ge (c/4)\big\{(s/n)^{1/2}+\varepsilon\big\}.
\end{align}
Combining these two inequalities, for $n\ge (10+2\log (1/c))^2$, 
we get 
\begin{align}
    R_{\text{TV}}^{\text{HDC}}(n,\epsilon,\Delta_s^{k-1},
    \bar{\bs \theta}_n) \ge 
    (c/4)\big\{(s/n)^{1/2}+\varepsilon\big\}
\end{align}
for every $k\ge 1$ and every $\varepsilon\in[0,1]$. The same 
argument can be used to show that all the inequalities in
Theorem~\ref{theorem@2} are valid for Huber's deterministic 
contamination model as well. Since the inclusions
$\calM_n^{\text{HDC}}(\varepsilon,\bs \theta^*)\subset \calM_n^{\text{OC}}(\varepsilon,\bs \theta^*)\cap \calM_n^{\text{PC}}
(\varepsilon,\bs \theta^*) \subset\calM_n^{\text{AC}}(\varepsilon, 
\bs \theta^*)$ hold true, we conclude that the lower bounds 
obtained for HC remain valid for all the other contamination 
models and are minimax optimal. 

The main tool in the proof of Theorem~\ref{theorem@2} is the 
following result \cite[Theorem 5.1]{CGR15}. There is a universal 
constant $c_1>0$ such that  for every $\varepsilon\in[0,1)$, 
\begin{align}
    \inf_{\bar{\bs \theta_n}}\sup_{\calM_n^{\rm HC}(\varepsilon,\Delta)}
    \bP\big(d(\bar{\bs \theta}_n,{\bs \theta}^*)\ge w_d(\varepsilon,\Delta)\big)\ge c_1,
\end{align}
where $w_d(\varepsilon,\Delta)$ is the modulus of continuity defined by
$w_d(\varepsilon,\Delta) = \sup \{d(\bs \theta,\bs \theta'): d_{\text{TV}}(\bs \theta,\bs \theta')
\le \varepsilon/(1-\varepsilon)\}$. Choosing $\bs \theta$ and $\bs \theta'$ to differ on only to coordinates, one can check that, for any $\varepsilon\le 1/2$, 
$w_{\rm TV}(\varepsilon,\Delta_s^{k-1}) \ge \varepsilon$, 
$w_{\rm H}(\varepsilon,\Delta_s^{k-1}) \ge  (\varepsilon/2)^{1/2}$ and 
$w_{\mathbb{L}^2}(\varepsilon,\Delta_s^{k-1}) \ge \sqrt{2}\varepsilon$. 
Combining with the lower bounds in the non-contaminated setting, this result
yields the claims of Theorem~\ref{theorem@2}. 
In addition, \eqref{WqTV} combined with the results of this section implies 
that the rate in \eqref{Wq} is minimax optimal.

\subsection{Confidence regions} 

We established so far bounds for the expected value of estimation error. 
The aim of this section is to present bounds on estimation error of the 
sample mean holding with high probability. This also leads to constructing confidence regions for the parameter vector $\btheta^*$. 
To this end, the contamination rate $\varepsilon$ and the sparsity $s$ 
are assumed to be known. It is an interesting open question whether 
one can construct optimally shrinking confidence regions for unknown
$\varepsilon$ and $s$.

\begin{thm}
    Let $\delta\in(0,1)$ be the tolerance level. If $\btheta^*\in \Delta_s^{k-1}$, then under any contamination model, the regions of $\Delta^{k-1}$
    defined by each of the following inequalities 
    \begin{align}
        &d_{\mathbb L^2}(\bar\bX_n,\bs\theta)\le (1/n)^{1/2} + \sqrt{2}\,\epsilon + \big(\log(1/\delta)/n\big)^{1/2},\\
        &d_{\textup{TV}}(\bar\bX_n,\bs\theta)\le (s/n)^{1/2} + 2\,\epsilon + \big(2\log(1/\delta)/n\big)^{1/2},\\
        &d_{\textup{H}}(\bar\bX_n,\bs\theta)\le 
        3.2 \big((s/n)\log(2s/\delta)\big)^{1/2} +(2\varepsilon)^{1/2} +  \big((1/n)\log(2/\delta)\big)^{1/2},
    \end{align}
    contain $\bs\theta^*$ with probability at least $1-\delta$.
\end{thm}

\begin{figure}
    \centering
    \includegraphics[width=.3\linewidth,valign=t]{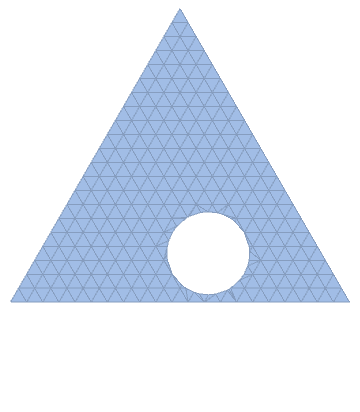}\quad
    \includegraphics[width=.3\linewidth,valign=t]{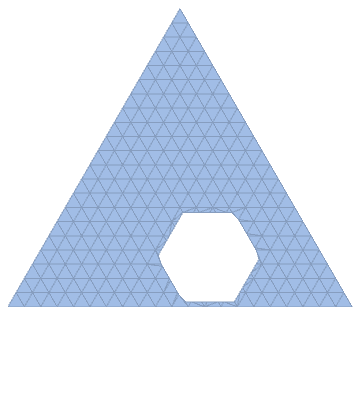}\quad
    \includegraphics[width=.3\linewidth,valign=t]{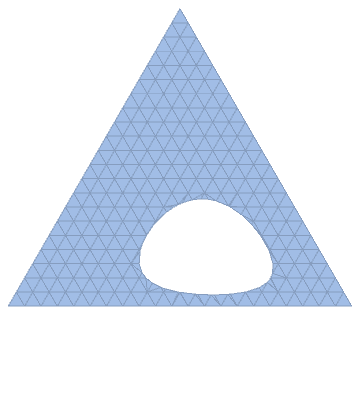}
    \vglue-20pt
    \caption{\label{fig:confidence} 
    The shape of confidence sets (white regions) for the distances $\mathbb L^2$ (left), TV (center), and Hellinger (right) when the sample mean is $(1/3,1/2,1/6)$.}
\end{figure}

To illustrate the shapes of these confidence regions, we 
depicted them in Figure~\ref{fig:confidence} for a three dimensional 
example, projected onto the plane containing the probability simplex. 
The sample mean in this example is equal to $(1/3,1/2,1/6)$. 

\section{Illustration on a numerical example}\label{sec:5}

We provide some numerical experiments which illustrate theoretical results of \Cref{sec:4}. The data set is the collection of 38 books written by Alexandre Dumas~(1802-1870) and 38 books written by Emile Zola~(1840-1902)\footnote{The works of both authors 
are available 
from \url{https://www.gutenberg.org/}}. To each author, we assign a parameter vector corresponding to the distribution of the number of words contained in the sentences used in the author's books. To be more clear, 
a sentence containing $l$ words is represented by vector $\be_l$, and if the parameter vector of an author is $(\theta_1,\dots,\theta_k)$, it means that 
a sentence used by the author is of size $l\in\{1,\dots,k\}$ with probability
$\theta_l$. We carried out synthetic experiments in which the reference parameter to estimate is the probability vector of Dumas, while the distribution of outliers is determined by the probability vector of Zola. Ground truths for 
these parameters are computed from the aforementioned large corpus of 
their works. Only the dense case $s=k$ were considered. For various values 
of $\varepsilon$ and $n$, a contaminated sample was 
generated by randomly choosing $n$ sentences either from Dumas' works 
(with probability $1-\epsilon$) or from Zola's works (with probability
$\epsilon$). The sample mean was computed for this corrupted sample, and 
the error with respect to Dumas' parameter vector was measured by the three  distances TV, $\mathbb L^2$ and Hellinger. This experiment was repeated 
$10^4$ times for each special setting to obtain information on error's
distribution. Furthermore, by grouping nearby outcomes we created samples of
different dimensions for illustrating the behavior of the error as a function 
of $k$.

\begin{figure}
\centering
\includegraphics[width=.45\textwidth]{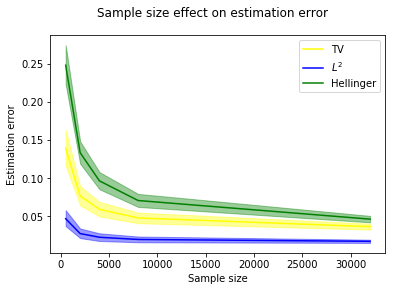}\quad
\includegraphics[width=.45\textwidth]{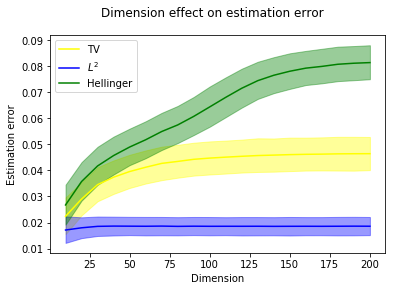}
\caption{\label{1} Estimation error of $\bar{\bs X}_n$ measured by total 
variation, Hellinger, and $\mathbb L^2$ distances as a function of (left panel) 
number of observations with contamination rate 0.2 and dimension $10^2$ 
and (right panel) dimension with contamination rate 0.2 and $10^4$ samples.
The interval between 5th and 95th quantiles of the error, obtained from $10^4$
repetitions, is also depicted for every graph. }
\end{figure}

\begin{figure}
\centering
\includegraphics[width=.45\linewidth,valign=t]{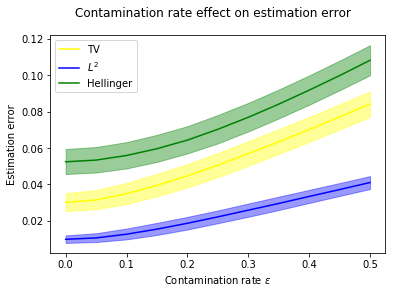}\quad
\includegraphics[width=.45\linewidth,valign=t]{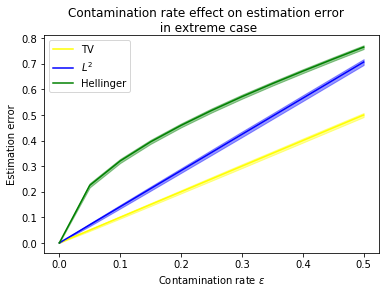}
\caption{\label{2} The estimation error of $\bar{\bs X}_n$ measured by total
variation, Hellinger, and $\mathbb L^2$ distances in terms of the contamination rate (with dimension $10^2$ and $10^4$ samples). At right, we plotted the error with respect to the contamination rate for an extreme case, where the reference and contamination distributions have the largest distance.
The interval between 5th and 95th quantiles of the error, obtained from $10^4$ trials, is also depicted.}
\end{figure}

The error of $\bar X_n$ as a function of the sample size~$n$, dimension~$k$, 
and contamination rate~$\epsilon$ is plotted in \Cref{1,2}. These plots are
conform to the theoretical results. Indeed, 
the first plot in \Cref{1} shows that the errors for the three distances is decreasing w.r.t.\  $n$. Furthermore, we see that up to some level of $n$ this decay is of order $n^{-1/2}$. The second plot in \Cref{1} confirms that the risk grows linearly in $k$ for 
the TV and Hellinger distances, while it is constant for the $\mathbb L^2$ error. 

Left panel of \Cref{2} suggests that the error grows linearly in terms of contamination rate. This is conform to our results for the TV and 
$\mathbb L^2$ errors. But it might seem that there is a disagreement with 
the result for the Hellinger distance, for which the risk is shown to 
increase at the rate $\varepsilon^{1/2}$ and not $\epsilon$. This is 
explained by the fact that the rate $\varepsilon^{1/2}$ corresponds to 
the worst-case risk,  whereas here, the setting under experiment does not
necessarily represent the worst case.  When the parameter vectors of the
reference and contamination distributions, respectively, are 
$\be_i$ and $\be_j$ with $i\neq j$ (\textit{i.e.}, when these two 
distributions are at the largest possible distance, which we call an 
extreme case), the graph of the error as a function of $\varepsilon$ 
(right panel of \Cref{2}) is similar to that of square-root function.

\section{Summary and conclusion}\label{sec:6}

We have analyzed the problem of robust estimation of the mean of a random vector
belonging to the probability simplex. Four measures of accuracy 
have been considered: total variation, Hellinger, Euclidean and Wasserstein  distances. In each case, we have established the minimax rates of the expected error of estimation under the sparsity assumption. In addition, confidence regions shrinking at the minimax rate have been proposed.  

An intriguing observation 
is that the choice of the distance has much stronger impact on the rate than 
the nature of contamination. Indeed, while the rates for the aforementioned distances are all different, the rate corresponding to one particular distance is not sensitive to the nature of outliers (ranging from Huber's contamination to the adversarial one). 
While the rate obtained for the TV-distance coincides with the previously known rate of robustly estimating a Gaussian mean, the rates we have established for the Hellinger and for the Euclidean distances appear to be new. Interestingly,
when the error is measured by the Euclidean distance, the quality of estimation
does not get deteriorated with increasing dimension. 
\tcbstoprecording

\appendix
\section{Proofs of propositions}\label{sec:7}


\begin{proof}[Proof of Proposition~\ref{prop@1} on page~\pageref{prop:1}]\label{proof:prop1}
Recall that $\hat O$ is the set of outliers in the Huber model. Let $O$ be any subset
of $\{1,\ldots,n\}$. It follows from the definition of Huber's model that if 
$\bP_n^O$ stands for the conditional distribution of $(\bX_1,\ldots,\bX_n)$ 
given $\hat O =O$, when $(\bX_1,\ldots,\bX_n)$ is drawn from $\bP^n\in\mathcal M_n^{\text{HC}}
(\varepsilon,\btheta^*)$, then $\bP_{n}^O\in\calM_n^{\text{HDC}} (|O|/n,\btheta^*)$. Therefore, 
for every $O$ of cardinality $o\ge 2\varepsilon n$, we have
\begin{align}
\sup_{\mathcal M_n^{\text{HC}}(\varepsilon,\btheta^*)}
\bfE[d(\hat\btheta_n,\btheta^*)\mathds 1(\hat O=O)]
		&=\sup_{\mathcal M_n^{\text{HC}}(\varepsilon,\btheta^*)}
			\bfE[d(\hat\btheta_n,\btheta^*)|\hat O=O]\bfP(\hat O=O)\\
		&\stackrel{}{\le}
		\sup_{\mathcal M_n^{\text{HDC}}(o/n,\btheta^*)}\bfE[d(\hat\btheta_n,\btheta^*)]
		\bfP(\hat O=O)\label{ineq:2}\\
		&\stackrel{(1)}{\le}\sup_{\mathcal M_n^{\text{HDC}}(1,\btheta^*)}\bfE[d(\hat\btheta_n
		,\btheta^*)]\bfP(\hat O=O).
\end{align}
Inequality (1) above is a direct consequence of the inclusion
$\mathcal M_n^{\text{HDC}}(o/n,\btheta^*)\subset \mathcal M_n^{\text{HDC}}(1,\btheta^*)$.
Summing the obtained inequality over all sets $O$ of cardinality 
$\ge 2\epsilon n$, we get
\begin{align}
\sup_{\mathcal M_n^{\text{HC}}(\varepsilon,\btheta^*)}
\bfE[d(\hat\btheta_n,\btheta^*)\mathds 1(|\hat O|\ge 2\varepsilon n)]
		&\le\sup_{\mathcal M_n^{\text{HDC}}(1,\btheta^*)}\bfE[d(\hat\btheta_n
		,\btheta^*)]\bfP(|\hat O|\ge 2\varepsilon n).
\end{align}
It follows from the multiplicative form of Chernoff's inequality 
that $\bfP(|\hat O|\ge 2\varepsilon n)\le e^{-n\varepsilon/3}$. 
This leads to the last term in inequality \eqref{ineq:1}.

Using the same argument as for \eqref{ineq:2}, for any $O$ of cardinality $o< 2n\varepsilon$,
we get
\begin{align}
\sup_{\mathcal M_n^{\text{HC}}(\varepsilon,\btheta^*)}
\bfE[d(\hat\btheta_n,\btheta^*)\mathds 1(|\hat O|<2n\varepsilon)]
		&\le 
		\sup_{\mathcal M_n^{\text{HDC}}(2\varepsilon,\btheta^*)}
		\bfE[d(\hat\btheta_n,\btheta^*)]
		\sum_{O:|O|\le 2n\varepsilon}\bfP(\hat O=O)\\
		&=\sup_{\mathcal M_n^{\text{HDC}}(2\varepsilon,\btheta^*)}
		\bfE[d(\hat\btheta_n,\btheta^*)].
\end{align}
This completes the proof of \eqref{ineq:1}. 

One can use the same arguments along with the Tchebychev inequality 
to establish \eqref{ineq:1b}. Indeed, for every $S$ of cardinality
$o\le 2\varepsilon n$, we have
\begin{align}
\sup_{\mathcal M_n^{\text{HC}}(\varepsilon,\btheta^*)}
r\bfP\big(d(\hat\btheta_n,\btheta^*)&>r \ \text{and}\ \hat O=O\big)\\
		&=\sup_{\mathcal M_n^{\text{HC}}(\varepsilon,\btheta^*)}
r\bfP\big(d(\hat\btheta_n,\btheta^*)>r \,|\,\hat O=O\big)
\bfP(\hat O=O)\\
		&\stackrel{}{\le}
		\sup_{\mathcal M_n^{\text{HDC}}(o/n,\btheta^*)}r
		\bfP\big(d(\hat\btheta_n,\btheta^*)>r\big)\bfP(\hat O=O)\\
		&\stackrel{}{\le}
		\sup_{\mathcal M_n^{\text{HDC}}(2\varepsilon,\btheta^*)}
		\bfE[d(\hat\btheta_n,\btheta^*)]\bfP(\hat O=O).
\end{align}
Summing the obtained inequality over all sets $O$ of cardinality 
$o\le 2\varepsilon n$, 
we get
\begin{align}
\sup_{\mathcal M_n^{\text{HC}}(\varepsilon,\btheta^*)}
r\bfP\big(d(\hat\btheta_n,\btheta^*)>r \ \text{and}\ |\hat O|\le 2\varepsilon n\big)
&\stackrel{}{\le}\sup_{\mathcal M_n^{\text{HDC}}(2\varepsilon,\btheta^*)}
		\bfE[d(\hat\btheta_n,\btheta^*)]\\
		& = R_d^{\text{HDC}}(n,2\varepsilon,\btheta^*,\hat\btheta).
\end{align}
On the other hand, it holds that
\begin{align}
\sup_{\mathcal M_n^{\text{HC}}(\varepsilon,\btheta^*)}
r\bfP\big(d(\hat\btheta_n,\btheta^*)&>r \ \text{and}\ |\hat O|> 2\varepsilon n\big)
\stackrel{}{\le} r\bfP\big(|\hat O|>2\varepsilon n) \le re^{-n\varepsilon/3},
\end{align}
and the claim of the proposition follows.
\end{proof}


\begin{proof}[Proof of Proposition~\ref{prop@2} on page~\pageref{prop:2}]\label{proof:prop2}
    Let $\btheta_1$ and $\btheta_2$ be two points in $\Theta$. We have 
\begin{align}
2R_d^{\text{HC}}(n,\varepsilon,\Theta,\hat\btheta_n)
		&\ge R_d^{\text{HC}}(n,\varepsilon,\btheta_1,\hat\btheta_n) + 
			R_d^{\text{HC}}(n,\varepsilon,\btheta_2,\hat\btheta_n)\\
		& \ge 
		\bfE_{(1-\varepsilon)\bP_{\btheta_1}+\varepsilon \bP_{\btheta_2}}[d(\hat\btheta_n,\btheta_1)] + 		
		\bfE_{(1-\varepsilon)\bP_{\btheta_1}+\varepsilon \bP_{\btheta_2}}[d(\hat\btheta_n,\btheta_2)].
\end{align}
To ease writing, assume that $n$ is an even number. Let $O$ 
be any fixed set of cardinality $n/2$. It is clear that the 
set of outliers $\hat O$ satisfies 
$$
p_O = \bP(\hat O = O) = \bP(\hat O = O^c) >0. 
$$
Furthermore, if $\bX_{1:n}$ is 
drawn from $((1-\varepsilon)\bP_{\btheta_1}+\varepsilon\bP_{\btheta_2})^{\otimes n}$, then 
its conditional distribution given $\hat O = O$ is exactly the same as the 
conditional distribution of $\bX_{1:n}\sim ((1-\varepsilon)\bP_{\btheta_2}+\varepsilon\bP_{\btheta_1})^{\otimes n}$ given $\hat O = O^c$. This implies that
\begin{align}
2R_d^{\text{HC}}&(n,\varepsilon,\Theta,\hat\btheta_n)\\
		& \ge p_O\big(\bfE_{(1-\varepsilon)\bP_{\btheta_1}+
		\varepsilon \bP_{\btheta_2}}
					[d(\hat\btheta_n,\btheta_1)|\hat O=O] + \bfE_{(1-\varepsilon)\bP_{\btheta_1}+
					\varepsilon \bP_{\btheta_2}}[d(\hat\btheta_n,\btheta_2)|\hat O =O^c]\big)\\
		& = p_O\bfE_{(1-\varepsilon)\bP_{\btheta_1}+\varepsilon \bP_{\btheta_2}}[d(\hat\btheta_n,\btheta_1)
				+d(\hat\btheta_n,\btheta_2)|\hat O=S] \ge  p_O d(\btheta_1,\btheta_2),
\end{align}
where in the last step we have used the triangle inequality. 
The obtained inequality being true for every $\btheta_1, \btheta_2 
\in \Theta$, we can take the supremum to get
\begin{align}
R_d^{\text{HC}}(n,\varepsilon,\Theta,\hat\btheta_n)		
		& \ge (p_O/2) \sup_{\btheta_1,\btheta_2\in\Theta} d(\btheta_1,\btheta_2) =+\infty.
\end{align}
This completes the proof.
\end{proof}

\section{Upper bounds on the minimax risk over the sparse simplex}

This section is devoted to the proof of the upper bounds on minimax risks 
in the discrete model with respect to various distances.

\begin{proof}[Proof of Theorem~\ref{theorem@1} on page~\pageref{thm:1}]\label{proof:thm1}
To ease notation, we set 
\begin{align}
    \bar{\bs X}_n=\frac1n\sum_i \bX_i,\quad 
    \bar{\bs Y}_n = \frac1n\sum_i \bY_i,\quad 
    \bar{\bs X}_{O} = \frac1o\sum_{i\in O} \bX_i,\quad  
    \bar{\bs Y}_{O} = \frac1o\sum_{i\in O} \bY_i
\end{align}
In the adversarial model, we have $\bX_i=\bY_i$ if $i\not\in O$
where $\bY_1,\dots,\bY_n$ are generated from the reference distribution $\bs \theta^*$. 
\begin{align}
d_{\mathbb L^2}(\bar{\bs X}_n,\bs \theta^*) &=
\big\|\bar{\bs X}_n - \bs \theta^*\big\|_{2} = 
\big\|\bar{\bs Y}_n - \bs \theta^* + \frac1n\sum_{i\in O} (\bX_i - \bY_i)\big\|_{2}\\
&\leq \big\|\bar{\bs Y}_n - \bs \theta^*\big\|_{2} + \frac{|O|}n \sup_{\bx,\by\in\Delta^{k-1}}\|\bx-\by\|_2\\
&=\big\|\bar{\bs Y}_n - \bs \theta^*\big\|_{2} + \sqrt{2}\varepsilon,
\end{align}
which gives us
\begin{align}
\sup_{\mathcal M_n^{\rm AC}(\varepsilon,\bs \theta^*)}
			\bfE[d_{\mathbb L^2}(\bar{\bs X}_n,\bs \theta^*)] &\leq \
			\sup_{\bs \theta^*}\bfE[d_{\mathbb L^2}(\bar{\bs Y}_n,\bs \theta^*)] + \sqrt{2}\varepsilon.
\end{align}
And for a fixed $\bs \theta^*$ it is well known that 
\begin{align}
    \bfE[d^2_{\mathbb L^2}(\bar{\bs Y}_n,\bs \theta^*)] 
    &= \sum_{j=1}^k \var[\bar{\bs Y}_{n,j}] = \sum_{j=1}^k \var\bigg[\frac1n\sum_{i=1}^n\fcar(Y_i=\bs e_j)\bigg] \\
    &= \frac1n\sum_{j=1}^k\var[\fcar(Y_1=\bs e_j)]\\
    &\le \frac1n\sum_{j=1}^k\bfE[\fcar(Y_1=\bs e_j)]= \frac1n.
\end{align}
Hence, we obtain $R_{\rm \mathbb{L}^2}^{\rm AC}(n,\epsilon,\Delta^{k-1}) \le (1/n)^{1/2} + \varepsilon$.
Similarly, 
\begin{align}
\dTV(\bar{\bs X}_n,\bs \theta^*) 
&\leq \big\|\bar{\bs Y}_n - \bs \theta^*\big\|_{1} + \frac{|O|}n \sup_{\bx,\by\in\Delta^{k-1}}\|\bx-\by\|_1\\
&=\big\|\bar{\bs Y}_n - \bs \theta^*\big\|_{1} + 2\varepsilon.
\end{align}
This gives 
\begin{align}
\sup_{\mathcal M_n^{\rm AC}(\varepsilon,\bs \theta^*)}
			\bfE[\dTV(\bar{\bs X}_n,\bs \theta^*)] &\leq \
			\sup_{\bs \theta^*}\bfE[\dTV(\bar{\bs Y}_n,\bs \theta^*)] + 2\varepsilon.
\end{align}
In addition, for every $\bs \theta^*$,
\begin{align}
    \bfE[\dTV(\bar{\bs Y}_n,\bs \theta^*)] &= 
    \frac{1}{2}\sum_{j=1}^{k}\bfE\big[\big|\bar{\bs Y}_{n,j}-\bs \theta^*_j\big|\big]\\
  &\leq \frac{1}{2}\sum_{j=1}^{k}\bigg(\bfE\big[\big|\bar{\bs Y}_{n,j}-\bs \theta^*_j\big|^2\big]\bigg)^{1/2}\\
  &= \frac{1}{2}\sum_{j=1}^{k}\bigg(\var\big[\bar{\bs Y}_{n,j}\big]\bigg)^{1/2}\\ &\stackrel{(1)}{=} \frac{1}{2}\sum_{j\in J} \big(\frac1n\bs\theta^*_j(1-\bs\theta^*_j) \big)^{1/2}\\
  &\stackrel{(2)}{\leq} \frac{1}{2}\,s^{1/2}\bigg(\sum_{j=1}^k \frac1n\bs \theta^*_j(1-\bs \theta^*_j)\bigg)^{1/2}\\ 
  &\leq \frac{1}{2}\big({s}/n\big)^{1/2},
\end{align}
where in (1) we have used the notation $J = \{j:\theta_j^*\neq 0\}$ and 
in (2) we have used the Cauchy-Schwarz inequality. This leads to 
\begin{align}
    R_{\rm TV}^{\rm AC}(n,\epsilon,\Delta^{k-1}) \le (k/n)^{1/2} + \epsilon. 
\end{align}
Finally, for the Hellinger distance
\begin{align}
d_{\rm H}(\bar{\bs X}_n,\bs \theta^*)&\le
d_{\rm H}(\bar{\bs X}_n,\bar{\bs Y}_n) + d_{\rm H}(\bar{\bs Y}_n,\bs \theta^*)\le
\sqrt2\dTV(\bar{\bs X}_n,\bar{\bs Y}_n)^{1/2} + d_{\rm H}(\bar{\bs Y}_n,\bs \theta^*),
\end{align}
where we have already seen that
\begin{align}
    \dTV(\bar{\bs X}_n,\bar{\bs Y}_n) &= \frac{o}{n}\big\|\bar{\bs X}_{O}-\bar{\bs Y}_{O}\big\|_{1} \le 2\epsilon.
\end{align}
This yields
\begin{align}
\sup_{\mathcal M_n^{\rm AC}(\varepsilon,\bs \theta^*)}
			\bfE[d_{\rm H}(\bar{\bs X}_n,\bs \theta^*)] &\leq \
			\sup_{\bs \theta^*}\bfE[d_{\rm H}(\bar{\bs Y}_n,\bs \theta^*)] + 2\sqrt{\varepsilon}.
\end{align}
Furthermore, for every $\bs\theta^*\in\Delta^{k-1}_s$,
\begin{align}
    \bfE[d^2_{\rm H}(\bar{\bs Y}_n,\bs \theta^*)] &= \bfE\bigg[\sum_{j=1}^{k}\Big(\sqrt{\bar{\bs Y}_{n,j}} - \sqrt{\bs \theta^*_j}\Big)^2\bigg] \\&\le
    \bfE\bigg[\sum_{j\in J}\frac{(\bar{\bs Y}_{n,j}-\bs \theta^*_j)^2}{\bs \theta^*_j}\bigg]\\
    &= \sum_{j\in J}\frac1{\bs \theta^*_j}\var\big[{\bar{\bs Y}_{n,j}}\big] = \sum_{j\in J}\frac1{\bs \theta^*_j}\frac{\bs \theta^*_j(1-\bs \theta^*_j)}{n} = \frac{s-1}{n}.
\end{align}
Hence, by Jensen's inequality $\bfE[d_{\rm H}(\bar{\bs Y}_n,\bs \theta^*)] < \sqrt{s/n}$. Therefore, we infer that
\begin{align}
R_{\rm H}^{\rm AC}(n,\epsilon,\Delta_s^{k-1}) \le (s/n)^{1/2} + 2\epsilon^{1/2}
\end{align}
and the last claim of the theorem follows. 
\end{proof}

\section{Lower bounds on the minimax risk over the sparse simplex}

This section is devoted to the proof of the lower bounds on minimax risks 
in the discrete model with respect to various distances. Note that the rates
over the high-dimensional ``sparse'' simplex $\Delta^{k-1}_s$
coincide with those for the dense simplex  $\Delta^{s-1}$. For this reason, 
all the lower bounds will be proved for $\Delta^{s-1}$ only (for $s\ge 2$
an even integer). 
In addition, we will restrict our attention to the distributions $\bs P$, 
$\bs Q$ over $\Delta^{s-1}$ that are supported by the set $\mathcal A$ 
of the elements of the canonical basis that is $\bs P(\mathcal A)=
\bs Q (\mathcal A)=1$. 

\begin{proof}[Proof of Theorem~\ref{theorem@2} on page~\pageref{thm:2}] \label{proof:thm2}
We denote by $\bs e_j$ the vector in $\mathbb R^s$ having all the 
coordinates equal to zero except the $j$th coordinate which is equal to one. 
Setting 
\begin{align}
    \bs\theta=\be_1,\quad\text{and} \quad
    \bs\theta'=\Big(1-\frac{\varepsilon}{1-\varepsilon}\Big)\be_1 + \frac{\varepsilon}{1-\varepsilon}\be_2    
\end{align}
we have 
\begin{align}
    \dTV(\bs \theta,\bs \theta') = \frac{\varepsilon}{1-\varepsilon},\quad
    d_{\mathbb L^2}(\bs \theta,\bs \theta') \ge \sqrt{2}\varepsilon,\quad 
    \text{and}\quad d_{\rm H}(\bs \theta,\bs \theta') \ge (\varepsilon/2)^{1/2}.
\end{align}
Therefore, modulus of continuity defined by 
$$
    w_d(\varepsilon,\Delta) = \sup \big\{d(\bs \theta,\bs \theta'): 
    \bs\theta,\bs\theta'\in\Delta,\ d_{\text{TV}}(\bs \theta,\bs \theta')\le \varepsilon/(1-\varepsilon)
    \big\}
$$ 
for a distance $d$ and a set $\Delta$, satisfies for any $\varepsilon\le 1/2$, 
\begin{align}\label{diff.modulus}
    w_{\rm TV}(\varepsilon,\Delta^{k-1}) \ge \varepsilon, \quad w_{\mathbb{L}^2}(\varepsilon,\Delta^{k-1}) \ge \sqrt{2}\varepsilon,
    \quad \text{and}\quad w_{\rm H}(\varepsilon,\Delta^{k-1}) 
    \ge  (\varepsilon/2)^{1/2}.
\end{align}
These bounds on the modulus of continuity are the first ingredient we need
to be able to apply Theorem 5.1 from \citep{CGR15} for getting  
minimax lower bounds. 

The second ingredient is the minimax rate in the non-contaminated case. 
are well known. For each of the considered distances, this rate is 
well-known. However, for the sake of completeness, we provide below a proof 
those lower bounds using Fano's method. For this, we use the Varshamov-Gilbert
lemma (see e.g. Lemma 2.9 in~\citep{T09}) and Theorem~2.5 in~\citep{T09}. 
The Varshamov-Gilbert lemma guarantees the existence of a set 
${\bs\omega}^{(1)}, \dots, {\bs\omega}^{(M)} \in \{0,1\}^{\lfloor s/2\rfloor}$ 
of cardinality~$M \geq 2^{s/16}$ such that 
\begin{align}
    \rho({\bs\omega}^{(i)},{\bs\omega}^{(j)}) \geq \frac{s}{16},\quad
    \text{for all}\quad i\neq j,    
\end{align} 
where $\rho(.,.)$ stands for the Hamming distance. Using these binary 
vectors $\{\bs \omega_j\}$, a parameter $\beta\in[0, \sqrt{n/s}]$ to be 
specified later and the ``baseline'' vector $\bs \theta^{(0)} 
= (1/s,\dots,1/s)$, we define hypotheses $\bs \theta^{(1)}$,
$\dots$, $\bs\theta^{(M)}$ by 
the relations
\begin{align}
  \bs \theta^{(i)}_{2j-1} = \bs \theta^{(0)}_{2j-1} + {\bs\omega}^{(i)}_j\frac{\beta}{\sqrt{ns}} \quad\text{ and }\quad \bs \theta^{(i)}_{2j} = \bs \theta^{(0)}_{2j}-{\bs\omega}^{(i)}_j\frac{\beta}{\sqrt{ns}} \quad \forall j \in \{0,\dots,\lfloor s/2\rfloor\}.  
\end{align}
Remark that~$\bs \theta^{(0)},\dots,\bs \theta^{(M)}$ are all probability 
vectors of dimension~$s$. Denoting the Kullback-Leibler divergence by 
$\dKL(.,.)$, one can check the conditions of Theorem~2.5 in~\cite{T09}:
\begin{align}
    d_{\mathbb L^2}(\bs \theta^{(i)}, \bs \theta^{(j)}) &\geq \sqrt2\frac{\beta}{\sqrt{ns}}\frac{\sqrt{s}}{4}=\frac{\beta}{\sqrt{8n}}
    \qquad \forall j\neq i,\\
    \dTV(\bs \theta^{(i)}, \bs \theta^{(j)}) &\geq \frac{\beta}{\sqrt{ns}}\frac{s}{16}=\frac{\beta\sqrt{s}}{16\sqrt{n}} \qquad \forall j\neq i,
\end{align}    
as well as
\begin{align}
    \dKL(\bs \theta^{(i)},\bs \theta^{(0)}) &\leq \sum_{j=1}^{\lfloor s/2\rfloor} \log\frac{\bs \theta^{(0)}_{2j-1}+\frac{\beta}{\sqrt{ns}}}{\bs \theta^{(0)}_{2j-1}}(\bs \theta^{(0)}_{2j-1}+\frac{\beta}{\sqrt{ns}}) 
    + \log\frac{\bs \theta^{(0)}_{2j}-\frac{\beta}{\sqrt{ns}}}{\bs \theta^{(0)}_{2j}}(\bs \theta^{(0)}_{2j}-\frac{\beta}{\sqrt{ns}})\\
    &\leq \sum_{j=1}^{\lfloor k/2\rfloor} \frac{\beta}{\sqrt{ns}\bs \theta^{(0)}_{2j-1}}(\bs \theta^{(0)}_{2j-1}+\frac{\beta}{\sqrt{ns}}) - \frac{\beta}{\sqrt{ns}\bs \theta^{(0)}_{2j}}(\bs \theta^{(0)}_{2j}-\frac{\beta}{\sqrt{ns}})\\
    &= \frac{\beta^2}{ns}\sum_{j=1}^{k}\frac1{\bs \theta^{(0)}_{j}} = \frac{\beta^2}{ns}\sum_{j=1}^{s}s = \frac{\beta^2s}{n} 
    \leq \frac{\alpha\log M}n \qquad \forall i\in\{1,\dots,M\},
  \end{align}
  for $\beta = \sqrt{\alpha}/4$. Now by applying the aforementioned theorem, 
  we obtain for the non-contaminated setting~($\varepsilon=0$)
\begin{align}
    \inf_{\bar{\bs \theta_n}}\sup_{\calM_n^{\rm HC}(0,\Delta^{s-1})}
    \bP\bigg(d_{\mathbb L^2}(\bar{\bs \theta}_n,{\bs \theta}^*)\ge \frac{\beta}{\sqrt{2n}}\bigg)
    \ge\frac{\sqrt{M}}{1+\sqrt{M}}\bigg(1-2\alpha-
    \sqrt{\frac{2\alpha}{\log M}}\bigg),\\
    \inf_{\bar{\bs \theta_n}}\sup_{\calM_n^{\rm HC}(0,\Delta^{s-1})}
    \bP\bigg(\dTV(\bar{\bs \theta}_n,{\bs \theta}^*)\ge \frac{\beta\sqrt{s}}{8\sqrt{n}}\bigg)\ge\frac{\sqrt{M}}{1+\sqrt{M}}\bigg(1-2\alpha-\sqrt{\frac{2\alpha}{\log M}}\bigg).
\end{align}
 Setting $M=2^{k/16}$ and $\alpha=1/32$, by Markov's inequality, one 
 concludes
 \begin{align}
 \inf_{\bar{\bs \theta}_n}R_{\mathbb{L}^2}^{\text{HC}}(n,0,\Delta^{s-1},\bar{\bs \theta}_n) &\ge c(1/n)^{1/2},\\
 \inf_{\bar{\bs \theta}_n}R_{\rm TV}^{\text{HC}}(n,0,\Delta^{s-1},\bar{\bs \theta}_n) &\ge c(s/n)^{1/2},
 \end{align}
 where $c=1/25600$. Since, $d_{\rm H}(\bs \theta,\bs \theta') \ge \dTV(\bs \theta,\bs \theta')$ for any distribution $\bs \theta$ and $\bs \theta'$
 \begin{align}
 \inf_{\bar{\bs \theta}_n}R_{\rm H}^{\text{HC}}(n,0,\Delta^{s-1},\bar{\bs \theta}_n) &\ge c(s/n)^{1/2}.
 \end{align}
 
 Finally, we apply \cite[Theorem 5.1]{CGR15} stating in our case for any distance $d$
 \begin{align}
 \inf_{\bar{\bs \theta}_n}R_{d}^{\text{HC}}(n,\varepsilon,\Delta^{s-1},\bar{\bs \theta}_n) &\ge c\big\{\inf_{\bar{\bs \theta}_n}R_{d}^{\text{HC}}(n,0,\Delta^{s-1},\bar{\bs \theta}_n) + w_{d}(\varepsilon,\Delta^{s-1})\big\}, 
 \end{align}
 for an universal constant $c$, which completes the proof of Theorem \ref{theorem@2}.
\end{proof}

\section{Proofs of bounds with high probability}

\begin{proof}
Suppose $\bX_i=\bY_i$ if $i\not\in O$ where $\bY_1,\dots,\bY_n$ are independently generated from the reference distribution $\bs P$ so that $\mathbf E[\bY_i] = \btheta^*$. For any $\bZ_1,\dots,\bZ_n$, let  $\bs\Phi_\square(\bZ_1,\dots,\bZ_n):= d_\square( \sum_{i=1}^n \bZ_i/n,\bs\theta^*)$, where $\square$ refers here to the distances $\mathbb L^2$ or TV. Given $\bY'_1,\dots,\bY'_n\in\Delta^{k-1}$ we have for every $i$
$$\bs\Phi_\square(\bY_1,\dots,\bY_i,\dots,\bY_n)-\bs\Phi_\square(\bY_1,\dots,\bY_{i-1},\bY'_i,\bY_{i+1},\dots,\bY_n) \leq \frac1n d_\square(\bY_i,\bY'_i).
$$
Furthermore, it can easily be shown that the last term  is bounded by $\sqrt{2}/n$ and $2/n$ for the distances $\mathbb L^2$ and TV, respectively. By bounded difference inequality (see for example Theorem~6.2 of~\cite{BLM13}) with probability at least $1-\delta$
\begin{align}
    \bs\Phi_{\mathbb L^2}(\bY_1,\dots,\bY_n) &\le \bfE\bs\Phi_{\mathbb L^2}(\bY_1,\dots,\bY_n) + \big(\log(1/\delta)/n\big)^{1/2}\\ 
    &\leq (1/n)^{1/2} + \big(\log(1/\delta)/n\big)^{1/2},\\
    \bs\Phi_{TV}(\bY_1,\dots,\bY_n)& \le \bfE\bs\Phi_{TV}(\bY_1,\dots,\bY_n) + \big(\log(2/\delta)/n\big)^{1/2}\\ 
    &\leq (s/n)^{1/2} + \big(\log(2/\delta)/n\big)^{1/2}.
\end{align}
Using $\bs\Phi_\square(\bX_1,\dots,\bX_n)\leq\bs\Phi_\square(\bY_1,\dots,\bY_n)+d_\square(\sum_{i\in O}\bX_i/n,\sum_{i\in O}\bY_i/n)$, one can conclude 
the proof of the first two claims of the theorem.

For the Hellinger distance, the computations are more tedious. We have to 
separate the case of small $\theta^*_j$. To this end, let $J = \{j:0<\theta^*_j<(1/n)\log(2s/\delta)\}$ and $J' = \{j:\theta^*_j\ge (1/n)\log(2s/\delta)\}$.
We have
\begin{align}
    \sum_{j\in J}\Big(\sqrt{\bar\bY_{n,j}} -\sqrt{\theta^*_j}\Big)^2 
    &\le \sum_{j\in J} (\bar\bY_{n,j} + \theta^*_j)\\
    &\le \frac1n\sum_{i=1}^n \Big(\sum_{j\in J} \bY_{i,j} - \theta_j^*\Big) + \frac{2s\log(2s/\delta)}{n}.
\end{align}
Since the random variables $U_i:=\Big(\sum_{j\in J} \bY_{i,j}\Big)$ are
iid, positive, bounded by $1$, the
Bernstein inequality implies that
\begin{align}
    \frac1n\sum_{i=1}^n \big( U_i - \mathbf E[U_1]\big)
    &\le   \sqrt{\frac{2\textbf{Var}(U_1)\log(2/\delta)}n} +
    \frac{\log(2/\delta)}{3n},
\end{align}
holds with probability at least $1-\delta/2$ for $0<\delta<1$.
One easily checks that $\sqrt{\textbf{Var}(U_1)}\le \sum_{j\in J}
\sqrt{\textbf{Var}(\bY_{1,j})}\le s \sqrt{\log(2s/\delta)/n}$. Therefore, with probability at least $1-\delta/2$, we have 
\begin{align}
    \frac1n\sum_{i=1}^n \Big(\sum_{j\in J} \bY_{i,j}-\theta_j^*\Big)
    &\le \frac{\sqrt{2}\,s\log(2s/\delta)}{n} + {\frac{\log(2/\delta)}{3n}}.
\end{align}
This yields
\begin{align}
    \sum_{j\in J} \Big(\sqrt{\bar \bY_{n,j}} -\sqrt{\theta^*_j}\Big)^2 
    &\le \frac{3.5 s{\log(2s/\delta)} + \log(2/\delta)}{n},\label{eq:43}
\end{align}
with probability at least $1-\delta/2$. 

On the other hand, we have
\begin{align}
    \sum_{j\in J'} \Big(\sqrt{\bar\bY_{n,j} } -\sqrt{\theta^*_j}\Big)^2 
        &\le \sum_{j\in J'} \frac{(\bar\bY_{n,j} -\theta^*_j)^2}{\theta_j^*}\\ 
        &\le s\max_{j\in J'} \frac{(\bar\bY_{n,j} -\theta^*_j)^2}{\theta_j^*}.
        \label{eq:44}
\end{align}
The Bernstein inequality and the union bound imply that, with probability 
at least $1-\delta/2$,
\begin{align}
    |\bar\bY_{n,j} -\theta^*_j|&\le \sqrt{\frac{2\textbf{Var}(\bY_{1,j})\log(2s/\delta)}{n}}
    +\frac{\log(2s/\delta)}{n},\qquad\forall j\in J'\\
    &\le\sqrt{\frac{2\theta_j^*\log(2s/\delta)}{n}}
    +\frac{\log(2s/\delta)}{n},\qquad\forall j\in J'\\
    &\le 2.5\sqrt{\frac{\theta_j^*\log(2s/\delta)}{n}},\qquad\forall j\in J'.
    \label{eq:45}
\end{align}
Combining \eqref{eq:44} and \eqref{eq:45}, we obtain
\begin{align}
    \sum_{j\in J'} \Big(\sqrt{\bar\bY_{n,j}} -\sqrt{\theta^*_j}\Big)^2 
        &\le 2.5^2\frac{s\log(2s/\delta)}{n},
        \label{eq:46}
\end{align}
with probability at least $1-\delta/2$. Finally, inequalities \eqref{eq:43} 
and \eqref{eq:46} together lead to  
\begin{align}
    d_{\textup{H}}^2(\bar{\bs Y}_n,\bs\theta^*)=\sum_{j=1}^n \Big(\sqrt{\bar\bY_{n,j}} -\sqrt{\theta^*_j}\Big)^2 
        &\le \frac{9.75 s\,\log(2s/\delta)}{n} + \frac{\log(2/\delta)}{n},
        \label{eq:46}
\end{align}
which is true with probability at least $1-\delta$. 
Using the triangle inequality and the fact that the Hellinger distance is
smaller than the square root of the TV-distance, we get
\begin{align}
    d_{\textup{H}}(\bar{\bs X}_n,\bs\theta^*)
    &\le d_{\textup{H}}(\bar{\bs X}_n,\bar{\bs Y}_n) + d_{\textup{H}}(\bar{\bs Y}_n,\bs\theta^*) \\
    &\le \sqrt{2\,d_{\textup{TV}}(\bar{\bs X}_n,\bar{\bs Y}_n)} + d_{\textup{H}}(\bar{\bs Y}_n,\bs\theta^*)\\
    &\le (2\varepsilon)^{1/2} +  3.2\sqrt{\frac{s\,\log(2s/\delta)}{n}} + \sqrt{\frac{\log(2/\delta)}{n}},
        \label{eq:47}
\end{align}
with probability at least $1-\delta$. This completes the proof of the theorem.
\end{proof}

\bibliography{refs}

\end{document}